\documentclass[3p,onecolumn]{article}
\usepackage{cite}
\usepackage{amsmath, amssymb}
\usepackage{amsfonts}
\usepackage{amsthm}
\usepackage[T1]{fontenc}
\usepackage[utf8]{inputenc} 
\newtheorem{theorem}{Theorem}[section]
\newtheorem{lemma}{Lemma}[section]

\newtheorem{problem}{Problem}
\newcommand{\mb}[1]{\mathbf{#1}}
\newcommand{\bs}[1]{\boldsymbol{#1}}
\newcommand{\cP}{\mathcal{P}}
\newcommand{\cC}{\mathcal{C}}
\newcommand{\cS}{\mathcal{S}}
\newcommand{\cJ}{\mathcal{J}}
\newcommand{\cO}{\mathcal{O}}

\newcommand{\cI}{\mathcal{I}}
\newcommand{\cE}{\mathcal{E}}
\newcommand{\bR}{\mathbb{R}}
\newcommand{\bC}{\mathbb{C}}
\newcommand{\os}{\overline{\mb{s}}}

\newcommand{\osigma}{\overline{\bs{\sigma}}}

\newcommand{\tp}{^{\mbox{\tiny\bf T}}}

\DeclareMathOperator*{\rk}{rk}

\begin{document}

\title{An identity theorem for the Fourier transform of polytopes on rationally parameterisable hypersurfaces}
\author{Konrad Engel \thanks{Universit\"at Rostock,  Institut f\"ur Mathematik, 18051 Rostock, Germany. E-mail: konrad.engel@uni-rostock.de}}

\maketitle

\begin{abstract}
A set $\cS$ of points in $\bR^n$ is called a rationally parameterisable hypersurface if $\cS=\{\bs{\sigma}(\mb{t}): \mb{t} \in D\}$, where $\bs{\sigma}: \bR^{n-1} \rightarrow \bR^n$ is a vector function with domain $D$ and rational functions as components. A generalized $n$-dimensional polytope in $\bR^n$ is a union of a finite number of  convex $n$-dimensional polytopes in $\bR^n$. The Fourier transform of such a generalized polytope $\cP$ in $\bR^n$ is defined by $F_{\cP}(\mb{s})=\int_{\cP} e^{-i\mb{s}\cdot\mb{x}} \,\mb{dx}$. We prove that $F_{\cP_1}(\bs{\sigma}(\mb{t})) = F_{\cP_2}(\bs{\sigma}(\mb{t}))\  \forall \mb{t} \in O$ implies $\cP_1=\cP_2$ if $O$ is an open subset of $D$ satisfying some well-defined conditions. Moreover we show that this theorem can be applied to quadric hypersurfaces that do not contain a line, but at least two points, i.e., in particular to spheres.
%
%
%
\end{abstract}

\section{Introduction}
\label{introduction}
A \emph{generalized $n$-dimensional polytope} in $\bR^n$ is a union of a finite number of  convex $n$-dimensional polytopes in $\bR^n$.
In this paper we study the \emph{Fourier transform} of $n$-dimensional generalized polytopes $\cP$ in $\bR^n$
\[
F_{\cP}(\mb{s})=\int_{\cP} e^{-i\mb{s}\cdot\mb{x}} \,\mb{dx}\,,
\]
restricted on rationally parameterisable hypersurfaces. Here bold symbols denote vectors in $\bR^n$ or $\bR^{n-1}$ and the product $\cdot$ is the standard scalar product.
A \emph{rationally parameterisable hypersurface} (briefly \emph{rp-hypersurface}) is a set $\cS$ of points in $\bR^n$ of the form
\begin{equation}
\label{rps}
\cS=\{\bs{\sigma}(\mb{t}): \mb{t} \in D\}\,,
\end{equation}
where 
\[
\bs{\sigma}(\mb{t})=
\begin{pmatrix}
\sigma_1(t_1,\dots,t_{n-1})\\
\vdots\\
\sigma_n(t_1,\dots,t_{n-1})
\end{pmatrix}\,,
\]
the functions
$\sigma_j$ are rational functions, $j=1,\dots,n$, and $D \subseteq \bR^{n-1}$ is the domain of $\cS$. In the case $n=2$ we use the notion \emph{rationally parameterisable curve} (briefly \emph{rp-curve}) instead of rp-hypersurface and use the notation $\cC$ instead of $\cS$.

For example, using polar resp. spherical coordinates and the standard substitution $t=\tan(\alpha)$, which implies $\cos(\alpha)=\frac{1-t^2}{1+t^2}$ and $\sin(\alpha)=\frac{2t}{1+t^2}$, one obtains that 
the unit circle 
\[
\cC_1=\left\{\left(\frac{1-t^2}{1+t^2}, \frac{2t}{1+t^2}\right): t \in \bR\right\}
\]
(with the missing point $(-1,0)$)
is an rp-curve and that the 3-dimensional unit sphere 
\[
\cS=\left\{\left(\frac{1-t_1^2}{1+t_1^2}\frac{1-t_2^2}{1+t_2^2}, \frac{2t_1}{1+t_1^2}\frac{1-t_2^2}{1+t_2^2}, \frac{2t_2}{1+t_2^2}\right): t_1 \in \bR, t_2 \in [-1,1]\right\}
\]
(with the missing segment $(-\sqrt{1-\lambda^2},0,\lambda)$, $\lambda \in (-1,1)$)
is an rp-hypersurface.
Clearly, also the hyperbola 
\[
\cC_2=\left\{\left(t, \frac{1}{t}\right): t \in \bR\setminus\{0\}\right\}
\]
and the parabola 
\[
\cC_3=\left\{(t, t^2): t \in \bR\right\}
\]
are rp-curves.
Since affine transformations do not violate the rationality this is also true for any circle, hyperbola, parabola and sphere. 

With the function $\bs{\sigma}: \bR^{n-1} \rightarrow \bR^n$ we associate its \emph{normalized function}, i.e., the function $\hat{\bs{\sigma}}: \bR^{n-1} \rightarrow \bR^{n-1}$ definded by

\[
\hat{\bs{\sigma}}(\mb{t})=
\begin{pmatrix}
\frac{\sigma_2(t_1,\dots,t_{n-1})}{\sigma_1(t_1,\dots,t_{n-1})}\\
\vdots\\
\frac{\sigma_n(t_1,\dots,t_{n-1})}{\sigma_1(t_1,\dots,t_{n-1})}
\end{pmatrix}\,.
\]

Note that $\hat{\bs{\sigma}}(\mb{t})$ is only defined if $\mb{t} \in D \setminus\sigma_1^{-1}(0)$, where $\sigma_1^{-1}(0)$ is the zero set of $\sigma_1$.
For a set $O \subseteq D\setminus\sigma_1^{-1}(0)$ let 
\[
\hat{\bs{\sigma}}(O)=\{\hat{\bs{\sigma}}(\mb{t}): \mb{t} \in O\}\,.
\]
In the following we need that, for an rp-hypersurface $\cS$ given by (\ref{rps}), an open subset $O$ of $D\setminus\sigma_1^{-1}(0)$ in $\bR^{n-1}$ satisfies two conditions:

\textbf{Hyperplane condition}: $\bs{\sigma}(O)$ is not contained in a hyperplane.

\textbf{Inner point condition}: There is a $\mb{t} \in O$ such that $\hat{\bs{\sigma}}(\mb{t})$ is an inner point of $\hat{\bs{\sigma}}(O)$ in $\bR^{n-1}$.

The main aim of the paper is the proof of the following theorem:
\begin{theorem}
\label{main_theorem}
Let $\cP_1$ and $\cP_2$ be generalized  $n$-dimensional polytopes in $\bR^n$ and let $\cS=\{\bs{\sigma}(\mb{t}): \mb{t} \in D\}$ be a rationally parameterisable hypersurface in $\bR^n$.
Assume that there exists an open subset $O$ of $D\setminus\sigma_1^{-1}(0)$ in $\bR^{n-1}$ that satisfies the hyperplane and the inner point condition.
If 
\[
F_{\cP_1}(\bs{\sigma}(\mb{t})) = F_{\cP_2}(\bs{\sigma}(\mb{t}))\quad \forall \mb{t} \in O\,,
\]
then $\cP_1=\cP_2$.
\end{theorem}

In particular we apply Theorem \ref{main_theorem} to quadric hypersurfaces, i.e., 
zero sets of an equation of the form
\begin{equation}
\label{hypersurface}
\frac{1}{2}\mb{s}\tp A \mb{s} + \mb{b}\tp \mb{s} + c = 0\,,
\end{equation}
where $A$ is a symmetric matrix different from the zero matrix.
\begin{theorem}
\label{quadric hypersurface theorem}
If $\cS$ is a quadric hypersurface that does not contain a line but at least two points, then, up to an exceptional set of  hypersurface measure zero, it is an rp-hypersurface with some parameter domain $D$ and every open subset $O$ of $D\setminus\sigma_1^{-1}(0)$ in $\bR^{n-1}$ satisfies the hyperplane and the inner point condition.
\end{theorem}

It is easy to construct a counterexample to the assertion in Theorem \ref{main_theorem} if $\bs{\sigma}(O)$ is contained in a hyperplane $h$ given by $\mb{a} \cdot \mb{x}=0$. If $\cP_2$ is the mirror of $\cP_1$ with respect to $h$, then, for all $\mb{s} \in h$, the product $\mb{s}\cdot \mb{x}$ is constant on the lines orthogonal to $h$, i.e., lines of the form 
$\{\mb{x}=\mb{x}_0+\lambda \mb{a}: \lambda \in \bR\}$, and hence by Cavalieri's principle
$F_{\cP_1}(\mb{s}) = F_{\cP_2}(\mb{s})$ for all $\mb{s} \in h$, and, in particular,
$F_{\cP_1}(\bs{\sigma}(\mb{t})) = F_{\cP_2}(\bs{\sigma}(\mb{t}))$ for all $\mb{t} \in O$.

In Section \ref{rational functions} we present several auxiliary results on rational functions and exponential functions.
The proof of Theorem \ref{main_theorem} is contained in Section \ref{proof of main theorem} and Theorem \ref{quadric hypersurface theorem} is proved in Section \ref{inner point condition}. Finally some open problems are presented in Section \ref{open problems}.

An important tool is the explicit representation of the Fourier transform of a convex polytope which is a central subject in the theory on the exponential valuation of polytopes built by Brion, Lawrence, Khovanskii,
Pukhlikov, and Barvinok \cite{Brion1988653, Lawrence1991259, Pukhlikov1992188, Barvinok1991149} and has origins in results of Motzkin and Schoenberg (mentioned by 
Davis \cite{Davis1964569}) as well of Grunsky \cite{Grunsky1955}. We recommend \cite{Barvinok2008, Beck2007} for studying this theory.

We show that the equality $F_{\cP_1}(\mb{s}) = F_{\cP_2}(\mb{s})$ for all $\mb{s} \in \bs{\sigma}(O)$ can be extended to $F_{\cP_1}(\mb{s}) = F_{\cP_2}(\mb{s})$ for all $\mb{s} \in \bR^n$. This is by no means clear in view of the preceding example. But having this extended equality the assertion in Theorem \ref{main_theorem} follows from the injectivity of the Fourier transform, see e.g. \cite{Plonka2018}.

The reconstruction of polytopes from sparse data of moments as well as from sparse data of the Fourier transform has been studied intensively in recent years. This concerns convex polytopes \cite{Milanfar1995, Gustafsson2000, Golub19991222, Elad20041814, Gravin2012596, Collowald20151079} and also non-convex polygons \cite{Plonka2013117} and generalized polytopes  \cite{Gravin2018}.

This paper should be considered a first step to solve Problem \ref{absolute_value}, i.e., the case that we do not know the values of the Fourier transform, but only the absolute value of it. In small-angle X-ray scattering and partly also in wide-angle X-ray scattering one can detect the absolute value of the Fourier transform on a sphere, the \emph{Ewald sphere}, see e.g. \cite{Raines2010, Barke2015}.
Therefore uniqueness questions and phase retrieval algorithms are important subjects for application in crystallography and other domains. Related results (without the restriction to hypersurfaces) can be found e.g. in \cite{Lemke2003, Beinert20151169, Beinert2018505}.

\section{Independence theorems over the field of rational functions}
\label{rational functions}
Let, as usual, $\bC[x]$ (resp. $\bC[\mb{x}]$) be the set of polynomials in the variable $x$ (resp. in the variables that are the components of $\mb{x}$) with complex coefficients. Moreover, let $\bC(x)$ (resp. $\bC(\mb{x})$) be the set of rational fractions in the variable $x$ (resp. in the variables that are the components of $\mb{x}$).
These polynomials and rational fractions define functions from $\bC$ (resp. $\bC^n$) into $\bC$, which we call also polynomials (resp. rational functions). For brevity we do not distinguish between the formal terms and the corresponding functions.

In the following let $\mathcal{O}$ be the zero function, i.e., the function that is everywhere 0 on a certain domain which is given by functions that are included in the concrete context.
In the same way we consider constant functions.
\begin{theorem}
\label{independence_theorem}
Let $p_k,r_k \in \bC(x), k \in [m]$, and let $r_k-r_{k'}$ be not constant for $k \neq k'$.
Let $I$ be an open interval in $\bR$.
If
\begin{equation}
\label{zero_sum}
\sum_{k=1}^m p_k(x)e^{r_k(x)}=0 \quad \forall x \in I\,,
\end{equation}
then 
\begin{equation}
\label{zero-function}
p_k =  \cO \quad \forall k \in[m]\,.
\end{equation}
\end{theorem}
\proof 
We may assume that $p_k \in \bC[x], k \in [m]$, i.e., these functions are polynomials, because, if necessary,  we can multiply (\ref{zero_sum}) by the main denominator of the rational functions $p_k , k \in [m]$.

Let $S$ be the (finite) set of all singularities of all involved rational functions $r_k$, $k \in [m]$.
Since all functions are analytic on $\bC \setminus S$, we have even
\begin{equation}
\label{zero_sum_complex}
\sum_{k=1}^m p_k(z)e^{r_k(z)}=0 \quad \forall z \in \bC \setminus S\,.
\end{equation}

Recall that we can write each rational function $r$ having $k$ singularities, i.e., poles including infinity, in the form
\begin{equation}
\label{partial_fractions}
r=c+\sum_{h=1}^g s^{(h)}\,,
\end{equation}
where $c$ is a constant and the items $s^{(h)}$ are the principal parts of the Laurent series about the pole $z_h$. Here each $s^{(h)}$, $h=1,\dots,g$, has the form
\[
s^{(h)}(z)=
\begin{cases}
\sum_{l=1}^{m_h} a_{-l}^{(h)}(z-z_h)^{-l}&\text{ if $z_h$ is a finite pole of order $m_h$}\,,\\
\sum_{l=1}^{m_h} a_{l}^{(h)}(z-z_h)^{l}&\text{ if $z_h=\infty$ is a pole of order $m_h$}\,.
\end{cases}
\]
The representation (\ref{partial_fractions}) can be obtained by decomposition of the rational function into partial fractions.
We allow \emph{dummy poles} of order 0, i.e., the corresponding principal part is the zero function, hence we may assume in the following that the functions $r_j$, $j=1,\dots,m$, have the same poles.

Let $\sigma_k$ be the sum of the orders of all poles of $r_k$ and let $\sigma=\sum_{k=1}^m \sigma_k$.
We proceed by induction on $\sigma$.
If $\sigma=0$, then $m=1$ and $r_1(x)$ is constant which implies that $p_1=\cO$. Thus we consider the induction step from $< \sigma$ to $\sigma$.

If $m=1$, then the assertion is trivially true. Thus let $m > 1$. 
We study an arbitrary pole $z_0$. In the following we assume that $z_0$ is finite. The case $z_0= \infty$ can be treated analogously or using the standard transformation $z \rightarrow \frac{1}{z}$.
Let $s_k$ be the principal part of the Laurent series of $r_k$, $k \in [m]$, about the pole $z_0$. 
Let $\omega$ be the maximal order of $z_0$ extended over all functions $r_k$, $k \in [m]$. Then we can write for all $k \in [m]$ the principal part $s_k$ in the form
\[
s_k(z)=\sum_{l=1}^{\omega} a_{-l,k} (z-z_0)^{-l}\,,
\]
where $a_{-\omega,k} \neq 0$ for at least one $k \in [m]$, say for $k=k^*$.




We put 
\[
z=z_0+\lambda e^{i \varphi}
\]
 with $\lambda, \varphi \in \bR$ and consider the limit process $\lambda \rightarrow 0$. 
Note that
\[
(z-z_0)^{-\omega}=\lambda^{-\omega} e^{-i\omega \varphi}\,.
\]
First we explain how we choose $\varphi \in [0, 2\pi)$. Let $k$ and $k'$ be different members of $[m]$ and let $l \in [\omega]$.
If $a_{-l,k} \neq a_{-l,k'}$, then there are only finitely many choices for $\varphi \in [0, 2\pi)$ such that
$\Re(a_{-l,k}e^{-il \varphi})= \Re(a_{-l,k'}e^{-i l\varphi})$ because $(a_{-l,k}-a_{-l,k'})e^{-il \varphi}$ has to be purely imaginary. 
Hence we can choose $\varphi \in [0, 2\pi)$ in such a way that
\begin{equation}
\label{real_part_injection}
 \Re(a_{-l,k}e^{-il \varphi}) =  \Re(a_{-l,k'}e^{-i l\varphi}) \Rightarrow
a_{-l,k} = a_{-l,k'} \quad \forall l \in [\omega]\,.
\end{equation}
Now we order the $\omega$-tuples
\[
\mb{a}_k=(\Re(a_{-\omega,k}e^{-i\omega \varphi}),\dots,(\Re(a_{-2,k}e^{-i 2 \varphi}),(\Re(a_{-1,k}e^{-i \varphi})), \quad k \in [m]\,,
\]
lexicographically. 
Obviously, there is some number $\kappa \in [m]$ such that $1 \le \kappa \le m$ and after an adequate renumbering
\[
\mb{a}_1=\dots = \mb{a}_{\kappa}
\]
and
\[
\mb{a}_l <_{lex} \mb{a}_1 \quad \forall l \in \{\kappa+1,\dots,m\}\,.
\]
In view of (\ref{real_part_injection}) we have
\[
s:=s_1 = \dots = s_{\kappa}\,.
\]
Let
\[
r_k'=r_k-s\,, \quad k \in [m]\,.
\]
Note that $r_k'$ is analytical at $z = z_0$ for all $k \in [\kappa]$.
Multiplying (\ref{zero_sum_complex}) by $e^{-s}$ leads to 
\begin{equation}
\label{r'}
\sum_{k=1}^m p_k(z)e^{r_k'(z)}=0 \quad \forall z \in \bC \setminus S\,.
\end{equation}
If 
\begin{equation}
\label{kappa zero}
\sum_{k=1}^{\kappa} p_k(z)e^{r_k'(z)}=0 \quad \forall z \in \bC \setminus S\,,
\end{equation}
then first $p_k =  \cO$ for all $k \in [\kappa]$ and then also  $p_k =  \cO$ for all $k \in \{\kappa+1,\dots,m\}$, i.e.,  (\ref{zero-function}), follow from the induction hypothesis. 
Thus assume that (\ref{kappa zero}) is not true. We can develop $\sum_{k=1}^m p_k(z)e^{r_k'(z)}$ into a power series having a first nonzero term $b_d (z-z_0)^d$, where $d \in\mathbb{N}$ and $b_d \neq 0$.
Let
\[
p_k'=(z-z_0)^{-d} p_k \,, \quad k \in [m]\,.
\]
From (\ref{r'}) we obtain
\begin{equation}
\label{p'}
\sum_{k=1}^m p_k'(z)e^{r_k'(z)}=0 \quad \forall z \in \bC \setminus S\,.
\end{equation}
In view of the power series expansion we have
\begin{equation}
\label{limit}
\lim \sum_{k=1}^{\kappa} p_k'(z)e^{r_k'(z)} =b_d \neq 0\,,
\end{equation}
where the limit can be arbitrarily taken for $z \rightarrow z_0$, but in particular also for $\lambda \rightarrow 0$.

Now let $k \in  \{\kappa+1,\dots,m\}$ and $\lambda \rightarrow 0$.

The principal part $s_k'$ of the Laurent series of $r_k'$ about $z=z_0$ has the form
\[
s_k'=\sum_{l=1}^{\omega_k} b_{-l,k} (z-z_0)^{-l}\,,
\]
where $1 \le \omega_k \le \omega$. Let  
\[
\beta_k=\Re(b_{-\omega_k,k}e^{-i \omega_k \varphi})\,.
\]
In view of the lexicographic ordering,
\[
\beta_k < 0\,.
\]
Obviously, for sufficiently small $\lambda > 0$,
\[
\Re(r_k'(z))=\beta_k \lambda^{- \omega_k} (1+O(\lambda)) \le \frac{\beta_k}{2} \lambda^{-1}\,.
\]
Consequently,
\[
|p_k'(z)e^{r_k'(z)}|=O(\lambda^{-d})|e^{\Re(r_k'(z))}|=O(\lambda^{-d})e^{\frac{\beta_k}{2} \lambda^{-1} } = o(1)\,.
\] 
Thus
\[
\sum_{k=\kappa+1}^{m} p_k'(z)e^{r_k'(z)}=o(1)\,.
\]
This is a contradiction to (\ref{p'}) and (\ref{limit}). \qed

The next lemma is folklore, but in order to make the paper self-contained we prove it.
\begin{lemma}
\label{polynomial_identity}
Let $P\in \bC(\mb{x})$. If there is an open subset $O$ of $\bR^n$ such that
\[
P(\mb{x})=0\quad \forall \mb{x} \in O\,,
\]
then $P=\cO$.
\end{lemma}
\proof 
First let $P$ be a polynomial, i.e., $P \in \bC[\mb{x}]$.
We  proceed by induction on $n$. The base case $n=1$ follows from the fundamental theorem of algebra. For the induction step $n-1 \rightarrow n$ we write $P$ in the form
\[
P=\sum_{j=0}^ka_j x_n^j\,,
\]
where $a_j\in \bC[x_1,\dots,x_{n-1}]$. Let $\pi$ be the projection from $\bR^n$ into $\bR^{n-1}$ that deletes the $n$-th component.
For each, but fixed $(x_1,\dots,x_{n-1}) \in \pi(O)$ there exist infinitely many values for $x_n$ such that
$(x_1,\dots,x_{n-1},x_n) \in O$. Hence, by the fundamental theorem of algebra
\[
a_j(x_1,\dots,x_{n-1})=0 \quad \forall (x_1,\dots,x_{n-1}) \in \pi(O) \text{ and } \forall j\,.
\]
Since $\pi(O)$ is an open subset in $\bR^{n-1}$, the induction hypothesis yields $a_j=\cO$ for all $j$ and hence also $P=\cO$.

The rational case $P \in \bC(\mb{x})$ can be reduced to the polynomial case by multiplication with the denominator.
\qed

\begin{theorem}
\label{rp-curve}
Let $\cS=\{\bs{\sigma}(\mb{t}): \mb{t} \in D\}$ be an rp-hypersurface and let $O$ be an open subset of $D$ in $\bR^{n-1}$ that  satisfies the hyperplane condition. Let $\mb{v}_k$, $k \in [m]$, be distinct points in $\bR^n$ and let $P_k \in \bC(\mb{x})$ for all $k \in[m]$.
If
\[
\sum_{k=1}^m P_k(\mb{x}) e^{-i \mb{v}_k \cdot \mb{x}} = 0 \quad \forall \mb{x} \in \bs{\sigma}(O)\,,
\]
then 
\[
P_k(\mb{x}) = 0 \quad \forall \mb{x} \in \bs{\sigma}(O)\,.
\]
\end{theorem}

\proof Let $\mb{x}_1$ be a fixed element of $\bs{\sigma}(O)$ and let $\mb{x}_1=\bs{\sigma}(\mb{t}_1)$. We have to show that $P_k(\mb{x}_1)=0$. By the supposition there are further points $\mb{x}_l \in \bs{\sigma}(O)$, $l=2,\dots,n+1,$ such that the points $\mb{x}_l, l \in [n+1]$, are not contained in a hyperplane.
Let $\mb{t}_l$, $l \in [n+1]$, be those points of $O$ for which $\mb{x}_l=\bs{\sigma}(\mb{t}_l)$. Let $\mb{t}_l=(t_{1,l},\dots,t_{n-1,l})\tp$. We may assume that the numbers $t_{1,l}, l \in [n+1]$, are distinct because otherwise we could disturb a little bit the numbers $t_{1,l}$ and would obtain the desired situation.
Let $q_j(t)$ be the interpolation polynomial such that, for all $l \in [n+1]$, we have $q_j(t_{1,l}) = t_{j,l}$, $j \in \{2,\dots,n-1\}$. Moreover, let $q_1(t)=t$. Further let $\mb{q}=(q_1,\dots,q_{n-1})\tp.$ 
Then $\mb{q}(t_{1,l})=\mb{t}_l$ for all $l \in [n+1]$.
Now we consider
\[
\bs{\rho}(t)=\bs{\sigma}(\mb{q}(t))\,.
\]
By construction,
\[
\bs{\rho}(t_{1,l}) = \mb{x}_l \quad \forall l \in [n+1]\,.
\]
Since $O$ is open there exists an open interval $I$ containing $t_{1,1}$ such that $\mb{q}(t) \in O$ for all $t \in I$.
Let for $k \in [m]$
\[
r_k=-i \mb{v}_k \cdot \bs{\rho} \text{ and } p_k=P_k(\bs{\rho})\,.
\]
By the assumption we have
\[
\sum_{k=1}^m p_k(t) e^{r_k(t)} = 0 \quad \forall t \in I\,.
\]
Of course, $p_k, r_k \in \bC(t)$. We further show that $r_k-r_{k'}$ is not constant for $k \neq k'$. Indeed, assume that $r_k-r_{k'}=c$. Then $(\mb{v}_k-\mb{v}_{k'}) \cdot \bs{\rho} = -\frac{c}{i}$ and hence the points $\mb{x}_l = \bs{\rho}(t_{1,l})$, $l \in [n+1]$, are contained in a hyperplane, a contradiction to the choice of these points. 
Finally we derive from Theorem \ref{independence_theorem} that $p_k=\cO$ and hence, in particular $0=p_k(t_{1,1})=P_k(\bs{\rho}(t_{1,1}))=P_k(\mb{x}_1)$.
\qed

A subset $S$ of $\bR^n$ is called a \emph{significant set} if 
$x_1 \neq 0$ for all $\mb{x} \in S$, and the associated set
\[
\hat{S}=\left\{\left(\frac{x_2}{x_1},\dots,\frac{x_n}{x_1}\right): \mb{x} \in S\right\}
\]
contains an open set in $\bR^{n-1}$. Note that $\bs{\sigma}(O)$ is a significant set if $O \subseteq D\setminus\sigma_1^{-1}(0)$ is open in $\bR^{n-1}$ and satisfies the inner point condition.

Recall that a function $f:\bR^n \rightarrow \bR$ is called \emph{homogeneous of degree $d$}, where $d$ is an integer, if $f(\lambda \mb{x}) = \lambda^d f(\mb{x})$ for all $\lambda \in \bR \setminus \{0\}$ and all $\mb{x}$ of the domain. We exclude here $\lambda=0$ because perhaps $f(\mb{0})$ is not defined.

\begin{lemma}
\label{homogeneous_lemma}
Let $P \in \bC(\mb{x})$ be homogeneous of degree $d$ and let $S$ be a significant subset of $\bR^n$.
If 
\begin{equation}
\label{hom0}
P(\mb{x})=0 \quad \forall \mb{x} \in S\,,
\end{equation}
then
\[
P = \cO\,.
\]
\end{lemma}
\proof 
Let $Q \in \bC(\mb{y})$ be defined by $Q(\mb{y})=P(1,y_1,\dots,y_{n-1})$, where $\mb{y}=(y_1,\dots,y_{n-1})$.
Because of $P(\mb{x})=x_1^d P(1,\frac{x_2}{x_1},\dots,\frac{x_n}{x_1})$ we have
\[
Q(\mb{y})=0 \quad \forall \mb{y} \in \hat{S}\,.
\]
Lemma \ref{polynomial_identity} implies $Q = \cO$ and thus also $P=\cO$. 
\qed

Let $\cE_d$ be the class of all function $F: \bC^n \rightarrow \bC$ of the form
\[
F(\mb{x}) = \sum_{k=1}^m P_k(\mb{x}) e^{-i \mb{v}_k \cdot \mb{x}} \,,
\]
where the points $\mb{v}_k$ are distinct, $k=1,\dots,m$, and the coefficients $P_k(\mb{x})$ are homogeneous rational functions of degree $d$.
We call these functions \emph{E-functions of degree $d$}. Note that any linear combination of E-functions of degree $d$ is again an E-function of degree $d$, i.e., these E-function form a vector space.

\begin{theorem}
\label{E-functions}
Let $\cS=\{\bs{\sigma}(\mb{t}): \mb{t} \in D\}$ be an rp-hypersurface and let $O$ be an open subset of $D$ in $\bR^{n-1}$ that  satisfies the hyperplane and the inner point condition.
If $F$ is an E-function and
\[
F(\mb{x}) = 0 \quad \forall \mb{x} \in \bs{\sigma}(O)\,,
\]
then 
\[
F=\cO.
\]
\end{theorem}
\proof The proof follows immediately from Theorem \ref{rp-curve} and Lemma \ref{homogeneous_lemma}. \qed

\section{Proof of Theorem \ref{main_theorem}}
\label{proof of main theorem}
First let $\cP$ be a convex $n$-dimensional polytope in $\bR^n$. Let $V_{\cP}$ be its vertex set. Here we consider each vertex as an element $\mb{v}$ of $\bR^n$.
Moreover let $L_{\cP}=\{\mb{v}-\mb{v}': \ \mb{v}, \mb{v}' \in V_{\cP}\}$.
From the theory on the exponential valuation of polytopes \cite{Brion1988653, Lawrence1991259, Pukhlikov1992188, Barvinok1991149} it is known that
\begin{equation*}
F_{\cP}(\mb{s}) = \sum_{\mb{v} \in V_{\cP}} Q_{\cP,\mb{v}}(\mb{s}) e^{-i \mb{v} \cdot \mb{s}} \quad \forall \mb{s} \in \bC^n \setminus Z_{\cP}\,,
\end{equation*}
where each $Q_{\cP,\mb{v}}(\mb{s})$ is a rational function of the form
\[
Q_{\cP,\mb{v}}(\mb{s})=\sum_{I \in \cI_{\cP}} \frac{\lambda_{\cP,\mb{v},I}}{\prod_{\bs{\ell} \in I} \bs{\ell} \cdot \mb{s}}\,,
\]
the numerators $\lambda_{\cP,\mb{v},I}$ are real numbers, $\cI_{\cP}$ is a family of $n$-element linearly independent subsets of 
$L_{\cP}$ and $Z_{\cP}$ contains those $\mb{s}$ for which a term $\bs{\ell} \cdot \mb{s}$ in the denominator is zero.
Note that $F_{\cP}$ is an E-function of degree $-n$.

Now let $\cP$ be a generalized $n$-dimensional polytope, i.e.,
\[
\cP=\cP_1\cup \dots \cup \cP_m
\]
with convex $n$-dimensional polytopes $\cP_1,\dots,\cP_m$ in $\bR^n$.
For $J \subseteq [m]$ let
\[
\cP_J=\bigcap_{j \in J} \cP_j
\]
and let $\cJ=\{J \subseteq [m]: \cP_J \text{ is $n$-dimensional}\}$.

Recall that the \emph{characteristic function} $\chi_T$ of a subset $T$ of $\bR^n$ is defined by
\[
\chi_T(\mb{x}) =
\begin{cases}
1&\text{if } \mb{x} \in T\,,\\
0&\text{otherwise}\,.
\end{cases}
\]
It is well-known that $1-\chi_{\cP_1\cup \dots \cup \cP_m} = \prod_{j=1}^m (1-\chi_{\cP_j})$ and
$\chi_{\cap_{j\in J}\cP_j}=\prod_{j\in J} \chi_{\cP_j}$
which leads to the inclusion-exclusion formula
\[
\chi_{\cP_1\cup \dots \cup \cP_m}=\sum_{\emptyset \neq J \subseteq[m]} (-1)^{|J|+1} \chi_{\cP_J}\,.
\]
Integration over $\bR^n$ and omission of zero terms yields
\[
F_{\cP}(\mb{s})=\sum_{J \in \cJ} (-1)^{|J|+1} F_{\cP_J}(\mb{s})\,,
\]
i.e., a linear combination of E-functions of degree $-n$, which is again an E-function of degree $-n$.




Now let $\cP_1$ and $\cP_2$ be the two generalized $n$-dimensional polytopes given in the theorem.

By the assumption of the theorem
\[
F_{\cP_1}(\bs{\sigma}(\mb{t})) - F_{\cP_2}(\bs{\sigma}(\mb{t}))=0\quad \forall \mb{t} \in O\,.
\]
Since $F_{\cP_1}-F_{\cP_2}$ is an E-function of degree $-n$, Theorem \ref{E-functions} implies
\[
F_{\cP_1}-F_{\cP_2}=\cO\,.
\]
Since the Fourier transform is obviously continuous in $\mb{s}$ we finally obtain
\[
F_{\cP_1}(\mb{s})-F_{\cP_2}(\mb{s}) =0 \quad \forall \mb{s} \in \bC^n\,.
\]
The injectivity of the Fourier transform, see e.g. \cite{Plonka2018}, implies $\cP_1=\cP_2$. \qed

\section{Proof of Theorem \ref{quadric hypersurface theorem}}
\label{inner point condition}
First we study the inner point condition more generally.
From the Implicit Function Theorem it immediately follows that the inner point condition is fulfilled for  the
open set $O \subseteq D\setminus\sigma_1^{-1}(0)$ if the determinant of the Jacobian is not everywhere 0, i.e.,
\begin{equation}
\label{Jacobian}
\exists \mb{t} \in O: \det\left(\frac{\partial \hat{\bs{\sigma}}(\mb{t})}{\partial \mb{t}}\right) \neq 0\,.
\end{equation}
But the explicit computation of this determinant is often difficult and thus we need some further sufficient conditions for (\ref{Jacobian}).
Let
\[
\os =\begin{pmatrix}
s_1\\\vdots\\s_{n-1}
\end{pmatrix}\text{  and  }
\osigma=
\begin{pmatrix}
\sigma_1\\\vdots\\\sigma_{n-1}
\end{pmatrix}\,.
\]
In the following we will assume that there is some $\mb{t}\in O$ such that
\begin{equation}
\label{first assumption}
\det\left(\frac{\partial \osigma(\mb{t})}{\partial \mb{t}}\right) \neq 0
\end{equation}
and discuss this assumption in detail at the end of this section.
Replacing $O$ by an adequate open subset, by the Implicit Function Theorem, we may even assume that (\ref{first assumption}) holds for all $\mb{t} \in O$ and that $\sigma$ maps $O$ bijectively onto
$\overline{O}=\osigma(O)$.
With the notation $f(\os)=\sigma_n(\osigma^{-1}(\os))$ we have for $\os=\osigma(\mb{t})$
\[
\sigma_n(\mb{t}) = f(\os)\,.
\]
Let finally
\[
\bs{\psi}(\os)=
\begin{pmatrix}
\frac{s_2}{s_1}\\
\vdots\\
\frac{s_{n-1}}{s_1}\\
\frac{f(s_1,\dots,s_{n-1})}{s_1}
\end{pmatrix}\,.
\]
Note that 
\[
\hat{\bs{\sigma}}(\mb{t}) = \bs{\psi}(\osigma(\mb{t}))\quad \forall \mb{t} \in O\,.
\]
Moreover,
\[
\det\left(\frac{\partial \hat{\bs{\sigma}}(\mb{t})}{\partial \mb{t}}\right) =
\det\left(\frac{\partial \bs{\psi}(\os)}{\partial \os}\right) 
\det\left(\frac{\partial \osigma(\mb{t})}{\partial \mb{t}}\right)\,.
\]
Thus, under the assumption (\ref{first assumption}) we only need that for $\os=\osigma(\mb{t})$
\begin{equation}
\label{second assumption}
\det\left(\frac{\partial \bs{\psi}(\os)}{\partial \os}\right) \neq 0\,.
\end{equation}
It is easy to check that
\[
\det\left(\frac{\partial \bs{\psi}(\os)}{\partial \os}\right) = \frac{1}{s_1^n}\left( \os \cdot \nabla f(\os) -f(\os)\right)\,.
\]
Let $I_{\os}=\{\lambda:  \lambda \os \in \overline{O}\}$ and let $\varphi: I_{\os} \rightarrow \bR$ be defined by
\[
\varphi(\lambda)=f(\lambda \os)\,.
\]
Then for $\os \in \overline{O}$ the following conditions are equivalent for $\lambda \in I_{\os}$:
\begin{align*}
\lambda\os \cdot \nabla f(\lambda\os) -f(\lambda\os)&=0\,,\\
\lambda \varphi'(\lambda)- \varphi(\lambda) &= 0\,,\\
\varphi(\lambda)&=c \lambda \text{  for some } c \in \bR\,,\\
f(\lambda \os)&=\lambda f(\os)\,.
\end{align*}
Hence we already proved the following theorem:
\begin{theorem}
\label{inner point theorem}
The inner point condition is fulfilled for  the
open set $O \subseteq D\setminus\sigma_1^{-1}(0)$ if there is some $\mb{t} \in O$ such that (\ref{first assumption}) holds
and with the point $\mb{s}=
\begin{pmatrix}
\os\\
 f(\os)
\end{pmatrix} \in \cS$
not all points $\lambda \mb{s}$ belong to $\cS$ if $\lambda$ runs through a sufficiently small neighborhood of 1.
\end{theorem}
\begin{lemma}
\label{rp-lemma}

Let $\cS = \{ \mb{s}\in \bR^n: \frac{1}{2}\mb{s}\tp A \mb{s} + \mb{b}\tp \mb{s} + c = 0\}$ be a quadric hypersurface that does not contain a line but at least two points.
Then $\rk(A|\mb{b})=n$ and  $\cS$ is an rp-hypersurface, possibly up to an exceptional set of hypersurface measure zero.
\end{lemma}
\proof
Assume that $\rk(A|\mb{b})<n$. Then there exists some $\mb{r}\neq \mb{0}$ in $\bR^n$ such that $A \mb{r} = \mb{0}$ and $\mb{b} \tp \mb{r} = 0$.
It follows that with $\mb{s} \in \cS$ the whole line $\{\mb{s}+\lambda \mb{r}: \lambda \in \bR\}$ belongs to $\cS$, a contradiction.

From linear algebra it is well-known that there exists an affine transformation $\mb{s}=T\mb{s}'+\mb{v}$ with a regular matrix $T$ and a shift vector $\mb{v}$ such that one of the following normal forms can be generated:

{\bf Case 1} $\rk(A)=n-1$.

\begin{equation}
\label{NF1}
\sum_{j=1}^{n-1} \epsilon_j s_j'^2=s_{n}'\,,\quad \text{ where } \forall j: \epsilon_j \in \{-1,1\}\,.
\end{equation}

{\bf Case 2} $\rk(A)=n$ and not all eigenvalues of $A$ have the same sign.

\begin{equation}
\label{NF2}
s_1's_n'+\sum_{j=2}^{n-1} \epsilon_j s_j'^2=c'\neq 0\,,\quad \text{ where } \forall j: \epsilon_j \in \{-1,1\}\,.
\end{equation}

Note that the line-condition implies $c'\neq0$ because otherwise $\cS$ would contain the line $s_2'=\dots=s_n'=0$.

{\bf Case 3} $\rk(A)=n$ and all eigenvalues of $A$ have the same sign.

\begin{equation}
\label{NF3}
\sum_{j=1}^{n} s_j'^2=1\,.
\end{equation}

In Cases 1 and 2 we can take the parameters
$t_j=s_j'$, $j \in[n-1]$, and thus $\mb{s}'$ is a function of $\mb{t}$ which we write in the form $\mb{s}'=\bs{\sigma}'(\mb{t})$. This leads to
$\bs{\sigma}(\mb{t}) = T \bs{\sigma}'(\mb{t}) + \mb{v}$.

In Case 3 we use spherical coordinates:
\begin{align*}
s_1'&=\cos(\varphi_1)\,,\\
s_2'&=\sin(\varphi_1)\cos(\varphi_2)\,,\\
s_3'&=\sin(\varphi_1)\sin(\varphi_2)\cos(\varphi_3)\,,\\
\vdots\\
 s_{n-1}'&=\sin(\varphi_1)\cdots\sin(\varphi_{n-2}) \cos(\varphi_{n-1})\,,\\
 s_{n}'&=\sin(\varphi_1)\cdots\sin(\varphi_{n-2}) \sin(\varphi_{n-1})\,,
\end{align*}
where $0 \le \varphi_j \le \pi$ for all $j \in [n-2]$ and $-\pi < \varphi_{n-1} \le \pi$.
Thus $\mb{s}'$ is a function of $\bs{\varphi}$ which we write in the form $\mb{s}'=\bs{\tau}(\bs{\varphi})$.
It is well-known that for
\[
t_j=\tan(\varphi_j/2) \text{ and } \varphi_j \neq \pi, \quad j \in [n-1]\,,
\]
i.e.,
\[
t_j\in \bR_{\ge 0}, \quad j \in [n-1],\,\text{ and } t_{n-1} \in \bR\,,
\]
\begin{align*}
\cos(\varphi_j)&=\frac{1-t_j^2}{1+t_j^2}\,,\\
\sin(\varphi_j)&=\frac{2 t_j}{1+t_j^2}\,,
\end{align*}
which shows that the components of $\mb{s}'$ and consequently also of $\mb{s}$ are rational functions
of $t_1,\dots,t_{n-1}$, up to the cases where one of the involved angles equals $\pi$. 

If we use the notation
\[
\bs{\arctan}(\mb{t})=\begin{pmatrix}
\arctan(t_1)\\
\vdots\\
\arctan(t_{n-1})
\end{pmatrix}\,
\]
and
\[
\bs{\sigma}'(\mb{t}) = \bs{\tau}(2\bs{\arctan}(\mb{t}))\,,
\]
then we have $\bs{\sigma}(\mb{t}) = T \bs{\sigma}'(\mb{t}) + \mb{v}$.

\qed

Note that we have $D=\bR^{n-1}$ in Cases 1 and 3 and $D=\bR_{\neq 0}\times \bR^{n-2}$ in Case 2.

We further restrict $\bs{\sigma}'$ to the first $n-1$ components and use the notation
\[
\os'=\osigma'(\mb{t}) = \begin{pmatrix}
\sigma_1'(\mb{t})\\
\vdots\\
\sigma_{n-1}'(\mb{t})
\end{pmatrix}\,.
\]

\begin{lemma}
\label{lemma hyp 1}
Let $\cS$ be a quadric hypersurface of the form of Lemma \ref{rp-lemma} and let $O \subseteq D\setminus\sigma_1^{-1}(0)$ be an open set. Then there is some $\mb{t} \in O$ such that
\begin{equation}
\label{determinant condition 2}
\det\left(\frac{\partial \osigma'(\mb{t})}{\partial \mb{t}}\right) \neq 0\,.
\end{equation}
\end{lemma}
\proof
We study the same cases as in the proof of Lemma \ref{rp-lemma}.
In Cases 1 and 2 the Jacobian $\frac{\partial \osigma'(\mb{t})}{\partial \mb{t}}$ is the identity matrix.
In Case 3 we have
\[
\frac{\partial \osigma'(\mb{t})}{\partial \mb{t}}=\frac{\partial \osigma'}{\partial \bs{\varphi}}
\frac{\partial \bs{\varphi}}{\partial \mb{t}}\,.
\]
Here the second matrix is obviously regular and for the first matrix it is easy to see that
\[
\det\left(\frac{\partial \osigma'}{\partial \bs{\varphi}}\right)=(-1)^{n-1} \prod_{j=1}^{n-1} \sin^{n-j}(\varphi_j)
\]
which cannot be identically zero if $\mb{t}$, i.e., also $\bs{\varphi}$, runs through an open set in $\bR^{n-1}$.
\qed

\begin{lemma}
\label{lemma hyp 2}
Let $\cS$ be a quadric hypersurface of the form of Lemma \ref{rp-lemma} and let $O \subseteq D\setminus\sigma_1^{-1}(0)$ be an open set.
If there is some $\mb{t} \in O$ such that
\begin{equation}
\label{determinant condition 3}
\det\left(\frac{\partial \osigma'(\mb{t})}{\partial \mb{t}}\right) \neq 0\,,
\end{equation}
then $O$ satisfies the hyperplane condition.
\end{lemma}
\proof The condition (\ref{determinant condition 2}) implies that $\osigma'(O)$ contains an open set. Now the normal forms (\ref{NF1})--(\ref{NF3}) imply that $\bs{\sigma}'(O)$ -- and hence also $\bs{\sigma}(O)$ -- cannot be contained in a hyperplane because the corresponding gradient is not constant. 
\qed


 
{\it Proof of Theorem \ref{quadric hypersurface theorem}}. 
By Lemmas \ref{lemma hyp 1} and Lemma \ref{lemma hyp 2} the open set $O$ satisfies the hyperplane condition. Moreover, 
with $\mb{s} \in \cS$
not all points $\lambda \mb{s}$ belong to $\cS$ if $\lambda$ runs through a sufficiently small neighborhood of 1. Indeed, if this would not be the case, then $\cS$ would contain the whole line $\{\lambda \mb{s}: \lambda \in \bR\}$ which is excluded.
In view of Theorem \ref{inner point theorem} it remains to check  (\ref{first assumption}).
Note that $\det\left(\frac{\partial \osigma(\mb{t})}{\partial \mb{t}}\right)$ is a rational function of $\mb{t}$. Lemma \ref{polynomial_identity} implies that if this determinant (multiplied by the main denominator)  would be zero for all $\mb{t} \in O$, then it would be zero for all $\mb{t} \in  D$.
Thus it is sufficient to prove that
\begin{equation}
\label{to show}
\exists \mb{t} \in D: \quad \det\left(\frac{\partial \osigma(\mb{t})}{\partial \mb{t}}\right) \neq 0\,.
\end{equation}
Note that
\begin{equation}
\label{chain rule}
\frac{\partial \osigma(\mb{t})}{\partial \mb{t}}=\frac{\partial \osigma(\mb{t})}{\partial \mb{s}'}\frac{\partial \bs{s}'}{\partial \mb{t}}=\overline{T} \frac{\partial \bs{s}'}{\partial \mb{t}}\,,
\end{equation}
where $\overline{T}$ is the minor of $T$ which can be obtained from $T$ by deleting the last row.
Let $\overline{\mb{T}}_j$ be its $j$-th column. 
Analogously, deleting the last component of $\mb{v}$ gives the vector $\overline{\mb{v}}$.

Further let 
\begin{align*}
\hat{\overline{T}}_j&=(\overline{\mb{T}}_1|\dots|\overline{\mb{T}}_{j-1}|\overline{\mb{T}}_n|\overline{\mb{T}}_{j+1}|\dots|\overline{\mb{T}}_{n-1})\,, \quad j \in [n-1]\,,\\
\hat{\overline{T}}_n&=(\overline{\mb{T}}_1|\dots|\overline{\mb{T}}_{n-1})\,.
\end{align*}

Since the transformation matrix $T$ is regular we have
$\rk(\overline{T})=n-1$
and hence there is some $j^* \in [n]$ such that
\begin{equation}
\label{rank=n-1}
\det(\hat{\overline{T}}_{j^*})\neq 0\,.
\end{equation}

The chain rule (\ref{chain rule}) implies

\begin{equation}
\label{Jacobian representation}
\frac{\partial \osigma(\mb{t})}{\partial \mb{t}} = \sum_{j=1}^{n}\overline{\mb{T}}_j \nabla(\bs{\sigma}_j'(\mb{t}))\tp\,.
\end{equation}
Now we distinguish between the cases in Lemma \ref{rp-lemma}.

{\bf Case 1 and 2}.

Since $\bs{\sigma}_j'(\mb{t})=t_j$ for $j \in [n-1]$
we have
\[
\frac{\partial \osigma(\mb{t})}{\partial \mb{t}}=\left(\mb{T}_1|\dots | \mb{T}_{n-1}\right)+\overline{\mb{T}}_n \nabla(\bs{\sigma}_n'(\mb{t}))\tp\,.
\]
In the following we  use (\ref{rank=n-1}) several times.

{\bf Case 1}.

Since $\bs{\sigma}_n'(\mb{t})=\sum_{j=1}^{n-1} \epsilon_j t_j^2$

\[
\nabla(\bs{\sigma}_n'(\mb{t})) = 2 \begin{pmatrix}\epsilon_1 t_1\\\vdots\\\epsilon_{n-1}t_{n-1} \end{pmatrix}
\]

which gives
\[
\frac{\partial \osigma(\mb{t})}{ \partial \mb{t}}=\left(\overline{\mb{T}}_j+2\epsilon_j t_j \overline{\mb{T}}_n\right)_{j \in[n-1]}\,.
\]
If $j^*=n$, then we take $\mb{t} = \mb{0}$ and
have at this point 
\[
\det\left(\frac{\partial \osigma(\mb{t})}{\partial \mb{t}}\right) = \det(\hat{\overline{T}}_{n}) \neq 0\,.
\]
If $j^*<n$, then we take $t_j=0$ for $j \in [n-1] \setminus \{j^*\}$ and $t_{j^*}= \lambda$ and have for these points
\[
\lim_{\lambda \rightarrow \infty} \frac{1}{\lambda}\det\left(\frac{\partial \osigma(\mb{t})}{\partial \mb{t}}\right) = 2\epsilon_{j^*}\det(\hat{\overline{T}}_{j^*}) \neq 0\,.
\]
Thus (\ref{to show}) is proved in this case.


{\bf Case 2}.

Since $\bs{\sigma}_n'(\mb{t})=\frac{1}{t_1}\left(c'-\sum_{j=2}^{n-1} \epsilon_j t_j^2\right)$ we have

\[
\nabla(\bs{\sigma}_n'(\mb{t})) = -\frac{1}{t_1} \begin{pmatrix} \sigma_n'(\mb{t})\\\epsilon_2 t_2\\\vdots\\\epsilon_{n-1}t_{n-1} \end{pmatrix}\,.
\]

We take $t_j=0$ for $j \in \{2,\dots,n-1\}$ and $t_{1}= \lambda$ and have for these points 
\[
\det\left(\frac{\partial \osigma(\mb{t})}{\partial \mb{t}}\right)=\det(\hat{\overline{T}}_{n}) - \frac{c'}{\lambda^2} \det(\hat{\overline{T}}_{1}) \,.
\]
If $j^*=n$, then we have
\[
\lim_{\lambda \rightarrow \infty}\det\left(\frac{\partial \osigma(\mb{t})}{\partial \mb{t}}\right)=\det(\hat{\overline{T}}_{n}) \neq 0\,,
\]
and if $j^*=1$, then we have
\[
\lim_{\lambda \rightarrow 0}\lambda^2\det\left(\frac{\partial \osigma(\mb{t})}{\partial \mb{t}}\right)=- c'\det(\hat{\overline{T}}_{1}) \neq 0\,.
\]
Thus we may assume that $\det(\hat{\overline{T}}_{1}) = \det(\hat{\overline{T}}_{n}) =0$ and that $1 < j^* < n$.
Here we take $t_{j^*}=\sqrt{\frac{c'}{\epsilon_j}}$ (possibly in $\bC$), $t_1=1$ and $t_j=0$ for $j \in [n]\setminus\{1,j^*\}$.
Then
\[
\det\left(\frac{\partial \osigma(\mb{t})}{\partial \mb{t}}\right)=\det(\hat{\overline{T}}_{n}) - \epsilon_{j^*}t_{j^*}\det(\hat{\overline{T}}_{j^*}) \neq 0\,.
\]
Thus (\ref{to show}) is proved also in this case.

{\bf Case 3}.

We have (with $\bs{\varphi} = 2\bs{\arctan}(\mb{t})$)
\[
\frac{\partial \mb{s}'}{\partial \mb{t}}=\frac{\partial \mb{s}'}{\partial \bs{\varphi}} \frac{\partial \bs{\varphi}}{\partial \mb{t}}
\]
and consequently by (\ref{chain rule})
\[
\det\left(\frac{\partial \osigma(\mb{t})}{\partial \mb{t}}\right)=\det\left(\overline{T}\frac{\partial \mb{s}'}{\partial \bs{\varphi}}\right) \det\left(\frac{\partial \bs{\varphi}}{\partial \mb{t}}\right)
\]
and thus it is sufficient to prove:
\begin{equation}
\label{to show case 3}
\exists \bs{\varphi} \in [0,\pi)^{n-2} \times (-\pi,\pi): \quad \det\left(\overline{T}\frac{\partial \mb{s}'}{\partial \bs{\varphi}}\right) \neq 0\,.
\end{equation}
With the notation
\[
\mb{a}(\varphi_1,\dots,\varphi_{n-1}) =
\begin{pmatrix}
-\sin(\varphi_1)\hfill\\
\cos(\varphi_1) \cos(\varphi_2)\hfill\\
\cos(\varphi_1) \sin(\varphi_2) \cos(\varphi_3)\hfill\\
\vdots\\
\cos(\varphi_1) \sin(\varphi_2) \sin(\varphi_3) \dots \sin(\varphi_{n-2}) \cos(\varphi_{n-1}) \\
\cos(\varphi_1) \sin(\varphi_2) \sin(\varphi_3) \dots \sin(\varphi_{n-2}) \sin(\varphi_{n-1})
\end{pmatrix}
\]
we have
\begin{multline*}
\frac{\partial \mb{s}'}{\partial \bs{\varphi}}=\\ \prod_{j=1}^{n-2} \sin^{n-j-1}(\varphi_j)
\begin{pmatrix}
\mb{a}(\varphi_1,\dots,\varphi_{n-1})&0&\dots&0&0\\
&\mb{a}(\varphi_2,\dots,\varphi_{n-1})&\dots&0&0\\
&&\dots&&\\
&&&\mb{a}(\varphi_{n-2},\varphi_{n-1})&0\\
&&&&\mb{a}(\varphi_{n-1})
\end{pmatrix}\,.
\end{multline*}
Note that
\[
\mb{a}\left(\frac{\pi}{2},\dots\right) =
\begin{pmatrix}
-1\\
0\\
\vdots\\
0\\
0
\end{pmatrix}\quad\text{and}\quad
\mb{a}\left(\frac{\pi}{4},\frac{\pi}{2},\dots,\frac{\pi}{2}\right) =
\begin{pmatrix}
-\frac{1}{2}\sqrt{2}\\
0\\
\vdots\\
0\\
\frac{1}{2}\sqrt{2}
\end{pmatrix}\,.
\]
If $j^*=n$, then we take $\varphi_1=\dots=\varphi_{n-1}=\frac{\pi}{2}$ and obtain for this point
\[
\det\left(\overline{T}\frac{\partial \mb{s}'}{\partial \bs{\varphi}}\right) = \det(-\hat{\overline{T}}_{n}) \neq 0\,.
\]
Thus we may assume that $1 < j^* < n$ and $\det(\hat{\overline{T}}_{n}) = 0$.
Now we take $\varphi_1=\dots=\varphi_{j^*-1}=\frac{\pi}{2}, \varphi_{j^*}=\frac{\pi}{4}, \varphi_{j^*+1}=\dots=\varphi_{n-1}=\frac{\pi}{2}$
and obtain
\[
\det\left(\overline{T}\frac{\partial \mb{s}'}{\partial \bs{\varphi}}\right) = \frac{1}{2} \sqrt{2} \left(\det(-\hat{\overline{T}}_{n})+\det(\hat{\overline{T}}_{j^*}) \right)\neq 0\,.
\]
Thus (\ref{to show}) is proved also in this final case.
\qed

\section{Open problems}
\label{open problems}
With the general definition of the Fourier transform of a function $f: \bR^n \rightarrow \bR$
\[
F_{f}(\mb{s})=\int_{\bR^n} f(\mb{x})e^{-i(\mb{s}\cdot\mb{x})} \,\mb{dx}
\]
Theorem \ref{main_theorem} is a result on characteristic functions $f=\chi_{\cP}$ of generalized $n$-dimensional polytopes $\cP$ in $\bR^n$. This result can be easily generalized to linear combinations of such functions, i.e., $f=\sum_{j=1}^m \lambda_i\chi_{\cP_j}$
because their Fourier transform is also an E-function.
\begin{problem}
Generalize Theorem \ref{main_theorem} to a larger class of functions $f: \bR^n \rightarrow \bR$.
\end{problem}
In Theorem \ref{main_theorem} we assume the equality $F_{\cP_1}(\bs{\sigma}(\mb{t})) = F_{\cP_2}(\bs{\sigma}(\mb{t}))$ for all $\mb{t} \in O$. By continuity, it is sufficient to require this equality only for the rational points of $O$, i.e., for a countable set of points.
\begin{problem}
Is there a finite set $T$ such that $F_{\cP_1}(\bs{\sigma}(\mb{t})) = F_{\cP_2}(\bs{\sigma}(\mb{t}))$ for all $\mb{t} \in T$ already implies $\cP_1=\cP_2$?
\end{problem}
It is easy to see that translations and point reflections of a polytope do not change the absolute value of its Fourier transform.
\begin{problem}
\label{absolute_value}
Assume that $|F_{\cP_1}(\bs{\sigma}(\mb{t}))| = |F_{\cP_2}(\bs{\sigma}(\mb{t}))|$ for all $\mb{t} \in O$. Under which conditions is it true that then $\cP_2$ can be obtained from $\cP_1$ by a translation and/or a point reflection?
\end{problem}
\section*{Acknowledgement}
I would like to thank Thomas Fennel and Stefan Scheel for presenting the motivation of this study and to Manfred Tasche for helpful hints and discussions. 

This work was partly supported by the European Social Fund (ESF) and the Ministry of Education, Science and Culture of Mecklenburg-Western Pomerania (Germany) within the project NEISS – Neural Extraction of Information, Structure and Symmetry in Images under grant no ESF/14-BM-A55-0006/19.

\bibliographystyle{plain}

\end{document}